\pdfoutput=1
\documentclass[12pt]{article}
\usepackage[utf8]{inputenc}
\usepackage[margin = 1in]{geometry}

\usepackage{graphicx}
\usepackage[all,cmtip]{xy}
\usepackage{subfigure,epsfig}
\usepackage[shortlabels]{enumitem}
\usepackage{amsmath}
\usepackage{amssymb}
\usepackage{amsthm}
\usepackage{bbm}
\usepackage[english]{babel}
\usepackage{fancyhdr}
\usepackage{stmaryrd}
\usepackage{tikz-cd}
\usepackage{mathtools}
\usepackage{titling}
\usepackage{dynkin-diagrams}
\usepackage[new]{old-arrows}
\usepackage{mathrsfs}
\usepackage[style=ext- alphabetic,sorting=nyt,maxnames=50,firstinits=true,maxalphanames=10,useprefix=true,articlein=false]{biblatex}

\usepackage[colorlinks=true, allcolors=blue]{hyperref}

\newtheorem{Theorem}{Theorem}[section]
\newtheorem{thm}{Theorem}[section]
\newtheorem{Corollary}[thm]{Corollary}
\newtheorem{Proposition}[thm]{Proposition}
\newtheorem{Lemma}[thm]{Lemma}
\newtheorem{Conjecture}[thm]{Conjecture}

\theoremstyle{definition}
\newtheorem{Definition}[thm]{Definition}

\theoremstyle{remark}
\newtheorem{Remark}[thm]{Remark}

\title{The Affine Closure of $T^*(\mathrm{SL}_{\lowercase{n}}/U)$}
\author{Boming Jia\thanks{\vspace{-0.1em}The author was supported by NSFC Grant No. 12225108 and the Shuimu Scholar Program in Tsinghua University.}}

\date{}

\newcommand{\Wedge}{\scalebox{0.8}{\raisebox{0.4ex}{$\bigwedge$}}}
\newcommand{\slc}{\mathfrak{sl}_3}

\newcommand{\slnb}{\mathfrak{sl}_{n-1}}

\newcommand{\so}{\mathfrak{so}}
\newcommand{\h}{\mathfrak{h}}

\newcommand{\Lieg}{\mathfrak{g}}

\newcommand{\Liet}{\mathfrak{t}}
\newcommand{\Lieh}{\mathfrak{h}}
\newcommand{\Lie}{\mathrm{Lie}}

\newcommand{\gl}{\mathfrak{gl}}
\newcommand{\slb}{\mathfrak{sl}_2}

\newcommand{\act}{\mathsf{act}}
\newcommand{\vol}{\mathsf{vol}}

\newcommand{\SL}{\mathrm{SL}}

\newcommand{\GL}{\mathrm{GL}}
\newcommand{\SO}{\mathrm{SO}}
\newcommand{\codim}{\mathrm{codim}}

\newcommand{\suchthat}{\,|\,}

\newcommand{\tr}{\mathrm{tr}}
\newcommand{\Tr}{\mathrm{Tr}}
\newcommand{\Id}{\mathrm{Id}}
\newcommand{\Sym}{\mathrm{Sym}}

\newcommand{\out}{\mathrm{out}}
\newcommand{\inn}{\mathrm{in}}

\DeclareMathOperator{\End}{End}
\DeclareMathOperator{\Hom}{Hom}
\DeclareMathOperator{\Spec}{Spec}

\DeclareMathOperator{\cspan}{span}

\newcommand{\Z}{\mathbb{Z}}
\newcommand{\C}{\mathbb{C}}

\newcommand{\Diff}{\mathcal{D}}

\newcommand{\Omin}{\mathcal{O}_{\textrm{min}}}
\newcommand{\Onilp}{\mathcal{O}_{\textrm{nilp}}}
\newcommand{\Ominbar}{\overline{\mathcal{O}}_\textrm{min}}
\newcommand{\Tsln}{\overline{T^*(\SL_n/U)}}
\newcommand{\Tslc}{T^*(\SL_3/U)}
\newcommand{\Tslcbar}{\overline{T^*(\SL_3/U)}}

\begin{document}

\setlength{\droptitle}{-6em}	

\maketitle

\vspace{-4em}

\begin{abstract}
We show that the affine closure $\overline{T^*(\mathrm{SL}_n/U)}$ has symplectic singularities, in the sense of Beauville. In the special case $n=3$, we show that the affine closure $\overline{T^*(\mathrm{SL}_3/U)}$ is isomorphic to the closure $\overline{\mathcal{O}}_\textrm{min}$ of the minimal nilpotent orbit $\mathcal{O}_{\textrm{min}}$ in $\mathfrak{so}_8$. Moreover, the quasi-classical Gelfand-Graev action of the Weyl group $W$ on $\overline{T^*(\mathrm{SL}_3/U)}$ can be identified with the restriction to $\overline{\mathcal{O}}_\textrm{min}$ of E.\ Cartan's triality action on $\mathfrak{so}_8$.   %<------------------- 
\end{abstract}

\section{Introduction}\

Let $G$ be a complex connected semisimple algebraic group, $B$ a Borel subgroup of $G$, and $U$ the unipotent radical of $B$. For example, in the case $G=\SL_n$, we can take $U$ to be the collection of all upper triangular matrices with all diagonal entries equal to $1$.
The homogeneous space $G/U$ is called the ``basic affine space". While $G/B$ is projective, the basic affine space $G/U$ is a quasi-affine variety. It turns out that many interesting problems in representation theory are related to the basic affine space. In particular, the algebra $\Diff(G/U)$ of algebraic differential operators on $G/U$ is well-studied, for example in  \cite{Gelfand}\cite{GK}.
In this paper, we study the total space of the cotangent bundle $T^*(G/U)$, of which the coordinate ring $\C[T^*(G/U)]$ is the quasi-classical counterpart of $\Diff(G/U)$. From a result by Ginzburg and Riche \cite{GR}, the coordinate ring $\C[T^*(G/U)]$ is finitely generated, and the affine closure of the basic affine space is defined as $$\overline{T^*(G/U)}\coloneqq \Spec \C[T^*(G/U)].$$
Symplectic singularities, a notion of which was first introduced by Beauville in \cite{B}, play an important role in representation theory. For instance, conic symplectic singularities (affine symplectic singularities with a good $\C^*$-action) admit universal flat Poisson deformations
and filtered quantizations.
There are many examples of symplectic singularities, for example, finite quotient singularities \cite{B}, normalization of the closure of nilpotent coadjoint orbits \cite{Pan}, and Nakajima quiver varieties \cite{BS}.
In \cite{GK}, Ginzburg and Kazhdan conjectured that
\begin{Conjecture}
    The affine closure $\overline{T^*(G/U)}$ has symplectic singularities.
\end{Conjecture}
And we prove the conjecture in section \ref{dq} for the special case $G=\mathrm{SL}_n$.

In the Lie algebra $\so_8$, there is a unique nilpotent adjoint orbit $\Omin\subset\so_8$ of minimal (positive) dimension $10$. The closure of the minimal orbit is $\Ominbar=\Omin\cup\{0\}.$ 
In section \ref{sl3}, we show that there is an isomorphism of affine varieties
    \begin{equation}\label{1.3}
        \overline{T^*(\SL_3/U)}\rightarrow\Ominbar.
    \end{equation}
In \cite{GK}, Ginzburg and Kazhdan constructed an action on $\overline{T^*(G/U)}$ by the Weyl group $W$ of $G$, called the ``Gelfand--Graev action". So in the case $G=\SL_3$ we have an $S_3$-action on $\Tslcbar$. On the other hand, the Lie algebra $\so_8$ has an $S_3$-symmetry called the triality action \cite{Cartan}\cite{MW}, and the restriction of the triality action gives an $S_3$-action on $\Ominbar$. In section \ref{triality}, we give a new interpretation of this triality action. In section \ref{ggaction}, we show that the isomorphism (\ref{1.3}) is $S_3$-equivariant.
In section \ref{symplectic}, we show the map (\ref{1.3}), when restricted on smooth points, is a symplectic isomorphism.

\bigskip
\noindent\textbf{Acknowledgments.}\ The author is grateful to both his doctorate advisor Victor Ginzburg and his postdoc advisor Peng Shan, for all their the encouragement, suggestions, and discussions. He would like to thank David Kazhdan, for the helpful correspondence containing statements related to the affine closure $\Tslcbar$, the minimal nilpotent orbit $\Omin$ in $\so_8$ and the triality action. He would also like to thank Ana Balibanu, Sam Evens, Yu Li, Daniil Rudenko, Minh-Tam Trinh and the referee for their useful suggestions and comments.

\section{The Affine Closure $\overline{T^*(\mathrm{SL}_{\lowercase{n}}/U)}$ Has Symplectic Singularities.}\label{dq}

We recall the following quiver theoretic construction of $\overline{T^*(\mathrm{SL}_{\lowercase{n}}/U)}$ in \cite{DKS}. A quiver $Q$ is a finite directed graph consisting of
a vertex set $I$ and an edge set $E$.
Write $Q^\textrm{op}=(I,\overline{E})$
for the opposite quiver obtained from $Q$ by reversing the orientation of edges.
The
double quiver of $Q$ is defined by $\overline Q=(I,E\sqcup \overline{E})$. 
A representation of $Q$ assigns a vector space $V(i)$ to each vertex $i\in I$ and a linear map $V(e):V(i)\rightarrow V(j)$ to each arrow $e\in E$ whose source and target are $i$ and $j$ respectively.
Let $Q$ be the following Dynkin quiver of type $\mathrm{A}_{n}$.
$$
\bullet \longrightarrow
	\bullet
	\longrightarrow
	\cdots
	\longrightarrow
	\bullet
	\longrightarrow
	\bullet
$$
Let $V=\bigoplus_{k=1}^{n-1}\Hom(\C^k,\C^{k+1})$, so each element of $V$ defines a representation of the quiver $Q$.
The cotangent space $T^*V$ is identified with $$\bigoplus_{k=1}^{n-1}\Hom(\C^k,\C^{k+1})\oplus\Hom(\C^{k+1},\C^{k})$$ via the  trace pairing
$$
    \sum_{k=1}^{n-1}\Tr (\beta_k\circ\alpha_k),
$$ for $\alpha_k\in\Hom(\C^k,\C^{k+1})$ and $ \beta_k\in\Hom(\C^{k+1},\C^{k})$.
So each element of $T^*V$ gives a representation of the double quiver $\overline{Q}$  
\begin{equation}\label{double quiver}
    \xymatrix@=3em{
{\C} \ar@<.5ex>[r]^{\alpha_1}
& {\C^2} \ar@<.5ex>[l]^{\beta_{1}} \ar@<.5ex>[r]^{\alpha_2}
& {\cdots} \ar@<.5ex>[l]^{\beta_{2}} \ar@<.5ex>[r]^{\alpha_{n-2}\ }
& {\C^{n-1}} \ar@<.5ex>[l]^{\beta_{n-2}\ } \ar@<.5ex>[r]^{\quad\alpha_{n-1}\ }
& {\C^{n}} \ar@<.5ex>[l]^{\quad\beta_{n-1}}}
\end{equation}
Throughout this paper, we use the expression $(\alpha,
\beta)$ to denote an element $$\bigoplus_{k=1}^{n-1}(\alpha_k,\beta_k)\in\bigoplus_{k=1}^{n-1}\Hom(\C^k,\C^{k+1})\oplus\Hom(\C^{k+1},\C^{k})=T^*V.$$
There is a natural action of $H\coloneqq\prod_{k=2}^{n-1}\SL_i(\C)$ on $V$ defined as follows.
Let $g=(g_2,\cdots,g_{n-1})\in H$, and $\alpha=(\alpha_1,\cdots,\alpha_{n-1})\in V$. Then we define 
\begin{equation}
g.\alpha=(g_{2}\circ\alpha_1,g_{3}\circ\alpha_2\circ g_{2}^{-1},\cdots,g_{n-1}\circ\alpha_{n-2}\circ g_{n-2}^{-1},\alpha_{n-1}\circ g_{n-1}^{-1}).
\end{equation}\label{eq:HonV}
This $H$-action on $V$
induces a Hamiltonian $H$-action on $T^*V$, of which the moment map $\mu:T^*V\longrightarrow\Lie(H)^*$ is given by
$$
    \mu_H(\alpha,\beta)(X)=\sum_{k=2}^{n-1}\Tr\left((\alpha_{k-1}\beta_{k-1}-\beta_k\alpha_k)X_k\right),
$$
where $(\alpha,\beta)=\bigoplus_{k=1}^{n-1}(\alpha_k,\beta_k)\in T^*V$, and $X=(X_2,X_3,\cdots,X_{n-1})\in\Lie(H)=\slb\times\slc\times\cdots\times\slnb$.

    We denote the zero fiber of the moment map by
    $$
        N\coloneqq\mu_H^{-1}(0)\subset T^*V.
    $$
    So $N$ is the subvariety consists of all the $(\alpha,\beta)$ in  $T^*V$ such that $\alpha_{k-1}\beta_{k-1}-\beta_k\alpha_k$ is a $k$-by-$k$ scalar matrix. Given $(\alpha,\beta)\in N$, and $k\in\{1,2,3,\cdots,n-1\}$, we define $\lambda_k$ as follows:
\begin{equation}\label{N}
    	\beta_k\alpha_k-\alpha_{k-1}\beta_{k-1}=\lambda_k\,\textrm{Id}_{\C^k}.
\end{equation}

\begin{Remark}
    A similar construction was applied in \cite{KP} to show the normality for the closure of every nilpotent orbit in $\gl_n$.
    Note that the scalars $\lambda_k$ depend on the choice of $(\alpha,\beta)\in N$, which is in contrast with the usual case of Nakajima quiver varieties. This difference is due to the definition of $H$ as the product of $\SL's$ rather than a product of $\GL's$ as usual. The choice of such $H$ also makes it much more difficult to construct an explicit resolution of singularities of $N\sslash H$ using Geometric Invariant Theory.
\end{Remark}

\begin{Theorem}[\textrm{Theorem 7.18} in \cite{DKS}]
The affine closure $\overline{T^*(\SL_n/U)}$ is isomorphic to the categorical quotient $N\sslash H$ as affine varieties.
\end{Theorem}

\begin{Lemma}
	Let $0=m_0\leq m_1\leq m_2\leq\cdots\leq m_{n-1}\leq m_n=n$.
	Suppose $m_k\leq k$ for all $k$. Then
	$$
	\sum_{k=1}^{n-1}m_k(m_{k+1}-m_k)\leq\sum_{k=1}^{n-1}k,
	$$
	and the equality holds if and only if $m_k=k$ for all $k$.
\end{Lemma}\label{ineq}
\begin{proof}
	Since $m_0=0$, we show that 
	\begin{align*}
		\sum_{k=0}^{n-1}m_k(m_{k+1}-m_k)
		=\sum_{k=0}^{n-1}\sum_{j=m_k+1}^{m_{k+1}}m_k
		\leq\sum_{k=0}^{n-1}\sum_{j=m_k+1}^{m_{k+1}}(j-1)
		=\sum_{k=1}^{n-1}k.\tag*{\qedhere}
	\end{align*}
\end{proof}
\begin{Remark}
    In fact, the LHS is the dimension of some (partial) flag variety of $\GL_n$, so it obtains maximum if and only if the flag variety is a complete flag variety, which has dimension equals to the RHS.
\end{Remark}

\begin{Lemma}\label{codim}
The singular locus of $\Tsln$ has codimension at least $4$.
\end{Lemma}
\begin{proof}
    Following \cite{DKS}, for each $\underline{m}=(m_1,\cdots,m_{n-1})$ satisfying the condition of Lemma \ref{ineq}, we set $$H(\underline{m})\coloneqq\prod_{k=1}^{n-1}\SL_{m_k}(\C),\ \ 
	\tilde{H}(\underline{m})\coloneqq\prod_{k=1}^{n-1}\GL_{m_k}(\C),\ \ V(\underline{m})\coloneqq\bigoplus_{k=1}^{n-1}\Hom(\C^{m_k},\C^{m_{k+1}}).$$
    Then we have an exact sequence
    $$
        1\rightarrow H(\underline{m})\rightarrow\tilde{H}(\underline{m})\xrightarrow{\varphi} T(\underline{m})\rightarrow 1,
    $$
    where $T(\underline{m})$ is a complex torus of rank $n-1$, and $\varphi$ is given by taking determinant of each product factor in $\tilde{H}(\underline{m})$.
     Let $\underline{n}=(1,2,3,\cdots,n-1)$.
	By Theorem 6.13 in \cite{DKS}, the affine variety $X$ can be written as disjoint union 
	   $$
	           X=\bigsqcup_{S,\delta}Q_{(S,\delta)}
	   $$
	where each $Q_{(S,\delta)}$ is a smooth hyperk\"ahler manifold that is a locally closed subset of $X$ labeled by an injective subrelation $S$ of $\leq$ on $\{1,2,\cdots,n-1\}$ and a function $\delta:\mathrm{dom}S\rightarrow\Z_{>0}$ which tells us which dimension vector $\underline{m}$ to use in the hyperkahler reduction construction of the strata $Q_{(S,\delta)}$, more precisely  $\underline{m}=\underline{n}-\sum_{(i,j)\in S}\delta(i)e_{ij},$ where $e_{ij}=(0,\cdots,0,1,\cdots,1,0,\cdots,0)$ is the $(n-1)$-tuple with $1$'s in positions in-between $i$ and $j$ and $0$'s otherwise.
 
	Let $x$ be a singular point of $\Tsln$. Then $x$ lies in some strata $Q_{(S,\delta)}$ other than the open dense generic strata $Q_{(\leq,0)}$. To the subrelation $S$, we associate the subtorus $T_S$ of $T(\underline{m})$ whose Lie algebra is $\Liet_S\coloneqq \cspan\{e_{ij}\mid i\leq_S j\}$.
 Now by Propostion 6.9 of \cite{DKS}, the strata $Q_{(S,\delta)}$ has dimension the same as the Hamiltonian reduction of $T^*V(\underline{m})$ by the subgroup $H_S\coloneqq\varphi^{-1}(T_S)\leq\tilde{H}(\underline{m})$. 
 Then $\dim H_S\geq\dim H(\underline{m})+1$.
	Thus we have
	\begin{align*}
		\dim Q_{(S,\delta)}& = 2(\dim V(\underline{m})-\dim H_S)\\
		& \leq 2(\dim V(\underline{m})-\dim H(\underline{m})-1)\\
		& =2\left(\sum_{k=1}^{n-1}m_km_{k+1}-(m^2_k-1)\right)-2\\
		& =2\left(n-2+\sum_{k=1}^{n-1}m_k(m_{k+1}-m_k)\right)
	\end{align*}
	Since $\underline{m}\neq(1,2,3,\cdots,n-1)$, by the previous lemma, we have $$\sum_{k=1}^{n-1}m_k(m_{k+1}-m_k)\leq\left(\sum_{k=1}^{n-1}k\right)-1.$$
	Thus
	$$
	\dim Q_{(S,\delta)}\leq2\left((n-3)+\sum_{k=1}^{n-1}k\right)
	=2\left(-1+\sum_{k=1}^{n}k\right)-4=\dim(\Tsln)-4.\qedhere
	$$
\end{proof}

\begin{Proposition}\label{symp}
	The smooth locus of\, $\Tsln$ admits a holomorphic symplectic form.
\end{Proposition}
\begin{proof}
	Let $X=\Tsln$. Let $\pi:N\rightarrow N\sslash H=X$ be the categorical quotient map. Let $X_\textrm{sm}$ be the smooth locus of $X$.
    We define the injective part $ N_\textrm{inj}$ of $N$ as 
    \begin{equation}\label{Nsurj}
        N_\textrm{inj} \coloneqq\{(\alpha,\beta)\in N\suchthat \textrm{all the $\alpha_k$ are injective}\}.
    \end{equation}
    Let $(\alpha,\beta)\in N_\textrm{inj}$. Then the stabilizer $H_\alpha$ for the $H$-action on $V$ defined in \eqref{eq:HonV} at $\alpha$ is trivial. So is the stabilizer $H_{(\alpha,\beta)}$ for the induced Hamiltonian $H$-action on $T^*V$ at $(\alpha,\beta)$. Now by Proposition 3.2 of \cite{HSS}, $(\alpha,\beta)$ is a smooth point of $N$.
    Moreover, by Theorem 4.5 of \cite{DKS}, the $H$-orbit of $\alpha$ in $V$ is closed. Now the $H$-orbit of $(\alpha,\beta)$ in $N$ is the graph of the morphism between the orbits $H.\alpha\rightarrow \beta.H$ given by $g\circ\alpha\mapsto \beta\circ g^{-1}$, hence also a closed $H$-orbit.
	So, by a Corollary of Luna's \'etale slice theorem (c.f. Proposition 5.7 in \cite{Dre}), we have $\pi(N_\textrm{inj})\subset X_\textrm{sm} \subset X$. From the proof of \textrm{Lemma 7.17} in \cite{DKS}, we know $N\setminus N_\textrm{inj}$ has codimension at least 2. So by the upper semicontinuity of the dimensions of the fibers of $\pi$ we have $$\codim(X_\textrm{sm}\setminus \pi(N_\textrm{inj}))\geq2.$$  
	Now by the proof of Proposition 7.2 in \cite{DKS}, we know $\pi(N_\textrm{inj})$ is identified with $$T^*(\mathrm{SL}_{\lowercase{n}}/U)\subset X_{sm}\subset X,$$ hence admitting a holomorphic symplectic form
	$\omega_0$. Then by Hartog's lemma, we can extend $\omega_0$ to a closed holomorphic 2-form $\omega$ on $X_\textrm{sm}$. By taking the top wedge of $\omega$, we have its points of degeneracy has codimension $1$, so $\omega$ has to be non-degenerated on $X_\textrm{sm}$, hence symplectic. 
\end{proof}

\begin{Definition}
        A normal variety $X$ is said to have symplectic singularities if
    \begin{enumerate}
        \item its smooth locus $X_\textrm{sm}$ carries a symplectic $2$-form $\omega$; and
        \item if $\nu:\widetilde{X}\rightarrow X$ is any resolution of singularities, then the pull-back $$\nu^*\omega\in\Omega^2(\nu^{-1}(X_\textrm{sm}))$$ extends to a holomorphic 2-form on $\widetilde{X}$.
    \end{enumerate}
\end{Definition}

\begin{Theorem}\label{main2}
	The affine closure $\overline{T^*(\SL_n/U)}$ has symplectic singularities.
\end{Theorem}
\begin{proof}
It is well-known that $X=\overline{T^*(\SL_n/U)}$ is normal (c.f. the remark after Definition 4.3 in \cite{Wang}.)
By Flenner's Theorem in \cite{Flenner}, it suffices to show that the smooth locus $X_\mathrm{sm}$ is a symplectic variety and $\codim(X\setminus X_\mathrm{sm})\geq4$.   
So the statement follows from Proposition \ref{symp} and Lemma \ref{codim}.
\end{proof}

\section{Isomorphism between $\overline{T^*(\SL_3/U)}$ and $\Ominbar\subset\so_8$}\label{sl3}

Let $m\geq4$.
Let $e_1,e_2,\cdots,e_{2m}$ be the natural basis of $\C^{2m}$.
Define an Euclidean inner product $(\ ,\ )$ on $\C^{2m}$ by
    \begin{equation}\label{inner}
        (e_i,e_j)=\delta_{i,2m+1-j}.
    \end{equation}
For $v_1\wedge v_2\in\Wedge^2\C^{2m}$, we define
\begin{align}\label{varphi}
	\qquad \varphi_{v_1\wedge v_2}:\ \C^{2m}\ &\longrightarrow \ \C^{2m}\\
	u\ \ &\longmapsto\ \, \left(v_1,u\right)v_2-
	\left(v_2,u\right)v_1.\nonumber
\end{align}
Extend by linearity, we get the following isomorphism of $\so_{2m}$-representations:
\begin{equation}\label{Lso}
    \Wedge^2\C^{2m}=\so_{2m}. 
\end{equation}

\begin{Definition}
    We say a subspace $W\subset \C^{2m}$ is isotropic if $(W,W)=0$. An element $f\in\Hom(\C^2,\C^{2m})$ is isotropic if its image is isotropic.
An element $\alpha\in\Wedge^2\C^{2m}$ is isotropic if for all $v_1^*, v_2^*\in (\C^{2m})^*$,
	$$
	    (\iota_{v_1^*}\alpha, \iota_{v_2^*}\alpha)=0.
	$$
\end{Definition}
\begin{Definition}
    We say $\alpha\in\Wedge^2\C^{2m}$ is decomposable if there exist some $v_1,v_2\in\C^{2m}$ such that $\alpha=v_1\wedge v_2$. 
\end{Definition}
\begin{Remark}
    By Pl\"ucker's Theorem we know that $\alpha\in\Wedge^2\C^{2m}$ is decomposable if and only if
$$\alpha\wedge\alpha=0\in\Wedge^4\C^{2m}.$$
\end{Remark}
\begin{Lemma}\label{iso}
	Let $v_1\wedge v_2\in\Wedge^2\C^{2m}$ 
	be a nonzero decomposable element. 
	Then 
	$$
	\mathrm{span} (v_1,v_2) \textrm{ is isotropic }
	\iff v_1\wedge v_2 \textrm{ is isotropic.}
	$$
\end{Lemma}
\begin{proof}
	($\Longrightarrow$)\ Clear.\\
	($\Longleftarrow$)\ It suffices to show that 
	$$
	(v_1,v_1)=(v_1,v_2)=(v_2,v_2)=0.
	$$
	Since $v_1\wedge v_2\neq0$, we have $v_1$ and $v_2$ are linearly independent. So there exist $v_1^*,v_2^*\in(\C^{2m})^*$ such that
	$$
	v_1^*(v_1)=v_2^*(v_2)=1,\textrm{ and } 
	v_1^*(v_2)=v_2^*(v_1)=0.
	$$
	Then
	\vspace{-1em}
	$$
	v_1=-\iota_{v_2^*}(v_1\wedge v_2)\textrm{, }
	v_2=\iota_{v_1^*}(v_1\wedge v_2),
	$$
	and our conclusion follows from the definition.
\end{proof}

\begin{Proposition}\label{ominiso} Let $\Ominbar$ be the  closure of the minimal nilpotent orbit $\Omin\subset\so_{2m}$. Then under the identification (\ref{Lso}),
	\begin{equation}\label{Omin}
	    \Ominbar=\{\alpha\in\Wedge^2\C^{2m}\suchthat \alpha\textrm{ is decomposable and isotropic}\}
	\end{equation}
\end{Proposition}
\begin{proof}
	Since $\Ominbar=\Omin\cup\{0\}$, it suffices to show that
	$$
	\Omin=\{\alpha\in\Wedge^2\C^{2m}\suchthat \alpha\textrm{ is a nonzero isotropic decomposable element}\}.
	$$
   By \textrm{Theorem 4.3.3} in \cite{CM}, the minimal orbit $\Omin$ is the adjoint orbit of the highest root vector $e_1\wedge e_2\in\Wedge^2\C^{2m}=\so_{2m}$.
	Since both decomposable and isotropic properties are invariant under the $\SO_{2m}$ action, and $e_1\wedge e_2$ is isotropic decomposable, so all elements in the minimal orbit are isotropic decomposable. But $\SO_{2m}$ acts transitively on the isotropic planes in $\C^{2m}$, and there is a scaling  $\C^*$-action on $\Omin$, so $\SO_{2m}$ acts transitively on the nonzero isotropic decomposable elements.
\end{proof}

Fix basis $e_1,e_2\in\C^2$, and define a symplectic bilinear form on $\Hom(\C^2,\C^{2m})$:
\begin{equation}\label{o1}
    \omega_1(f_1,f_2)=(f_1(e_1),f_2(e_2))-(f_1(e_2),f_2(e_1)).
\end{equation}
The natural $\SL_2$-action on $\Hom(\C^2,\C^{2m})$ is Hamiltonian with moment map
\begin{align*}
    & \mu_{\SL_2}:\Hom(\C^2,\C^{2m}) \longrightarrow \slb^*\cong\slb,\\
    & \mu_{\SL_2}(f) =
    \renewcommand*{\arraystretch}{1.5}
    \begin{pmatrix*}
    (f(e_1),f(e_2)) && \hspace{-1em}(f(e_1),f(e_1))\\
    (f(e_2),f(e_2)) && \hspace{-1em}-(f(e_1),f(e_2))
    \end{pmatrix*}.
\end{align*}
So the zero fiber of the $\SL_2$-moment map on \mbox{$\Hom(\C^2,\C^{2m})$ is}
    \begin{equation}
        N_1\coloneqq\mu_{\SL_2}^{-1}(0)=\{f\in\Hom(\C^2,\C^{2m})\suchthat f \textrm{ is isotropic with respect to }(\ ,\ )\}.
    \end{equation}
The following result is probably known to experts, but we still include a proof here for completeness. 
\begin{Proposition}\label{ominhamil} The Hamiltonian reduction $N_1\sslash\SL_2$ is isomorphic to the closure $\Ominbar$ of the minimal orbit in $\so_{2m}$ as affine varieties.
\end{Proposition}
\begin{proof}
Recall from Proposition \ref{ominiso}, $$\Ominbar=\{\alpha\in\Wedge^2\C^{2m}\suchthat \alpha\textrm{ is decomposable and isotropic}\}.$$
Denote
	\vspace{-1em}
	\begin{align*}
		& \left(\Wedge^2\C^{2m}\right)_\textrm{dec}\coloneqq\{\alpha\in\Wedge^2\C^{2m}\suchthat \alpha\mathrm{\ is\ decomposable}\},\\
		& \left(\Wedge^2\C^{2m}\right)_{\textrm{dec},\,\textrm{iso}}\coloneqq\{\alpha\in\Wedge^2\C^{2m}\suchthat \alpha\mathrm{\ is\ decomposable\ and\ isotropic}\}.
	\end{align*}
	By the First Fundamental Theorem of invariant theory for $\SL_2$ (\textrm{Theorem 2} in page 387 of \cite{Procesi}), we have $$\C[(\C^2)^*\otimes\C^{2m}]^{\SL_2}=\C\left[\left(\Wedge^2\C^{2m}\right)_\textrm{dec}\right].$$ 
	Since $\SL_2$ is reductive and $N_1\subset\Hom(\C^2,\C^{2m})=(\C^2)^*\otimes\C^{2m}$ is an $\SL_2$-invariant sub-variety, the restriction from $(\C^2)^*\otimes\C^{2m}$ to $N_1$ gives a surjective map
	\begin{align}\label{surj}
		\C\left[\left(\Wedge^2\C^{2m}\right)_\textrm{dec}\right]\quad& \longrightarrow\qquad\qquad\C[N_1]^{\SL_2}\\
		f\qquad \qquad& \longmapsto \quad\left(e_1^*\otimes v_1+e_2^*\otimes v_2\mapsto f(v_1\wedge v_2)\right)\nonumber
	\end{align}
By Lemma \ref{iso} we have the kernel of (\ref{surj}) is precisely the ideal $\mathcal{I}_\textrm{iso}$ of functions on $\left(\Wedge^2\C^{2m}\right)_\textrm{dec}$ vanishing identically on the isotropic elements. Hence
$$
    \C[N_1]^{\SL_2}=\C\left[\left(\Wedge^2\C^{2m}\right)_\textrm{dec}\right]/\mathcal{I}_\textrm{iso}=\C\left[\left(\Wedge^2\C^{2m}\right)_\textrm{dec,\,iso}\right].\qedhere
$$
\end{proof}

    Let $\eta:\C^3\oplus\C\oplus\C^*\oplus(\C^3)^*\longrightarrow\C^8$ be the isomorphism given by
    \begin{equation}\label{eta}
          \eta(z_1,z_2,z_3,a,b,w_1,w_2,w_3)=z_1e_2+z_2e_3-z_3e_8+ae_4
        +be_5-w_3e_1+w_2e_6+w_1e_7
    \end{equation}
    Then $\eta$ is an isometry with respect to the inner product on $\C^4\oplus(\C^4)^*$ given by the natural pairing of $\C^4$ and $(\C^4)^*$ and the inner product $(\,,\,)$ on $\C^8$ defined by taking $m=4$ in (\ref{inner}).
    We define a linear isomorphism
\begin{align*}
	F:
	T^*V& \longrightarrow \Hom(\C^2,\C^8)=\Hom(\C^2,\C^4\oplus(\C^4)^*)\\
	(\alpha,\beta) &\longmapsto\left(v\mapsto (\alpha_2\oplus-\beta_1)(v)\oplus\frac{(\beta_2\oplus\alpha_1)\wedge v}{e_1\wedge e_2}\right)
\end{align*}

\begin{Proposition}\label{Fthm}
	The map $F$ is an $\SL_2$-equivariant symplectic isomorphism between $T^*V$ and $(\Hom(\C^2,\C^8),\omega_1)$, where $\omega_1$ is given by (\ref{o1}).
\end{Proposition}
\begin{proof}
	First, we define $F_0:T^*V\longrightarrow T^*(\Hom(\C^2,\C^4))$ which maps $(\alpha,\beta)$ to the following element in $\Hom(\C^2,\C^4)\oplus\Hom(\C^4,\C^2)$:
 \begin{equation*}
    \xymatrix@=4em{
{\C^2} \ar@<.5ex>[r]^{\hspace{-2.5em}\alpha_2\oplus-\beta_1}
& {(\C^3\oplus\C)=\C^4} \ar@<.5ex>[l]^{\hspace{-2.5em}\beta_2\oplus\alpha_1}}
\end{equation*}
	Then $F_0$ is an $\SL_2$-equivariant symplectic isomorphism.
	Next, under the identifications
	\begin{align*}
		& T^*(\Hom(\C^2,\C^4))=T^*(\C^4\oplus\C^4)=(\C^4\oplus\C^4)\oplus((\C^4)^*\oplus(\C^4)^*),\\
		& \Hom(\C^2,\C^4\oplus(\C^4)^*)=(\C^4\oplus(\C^4)^*)\oplus(\C^4\oplus(\C^4)^*),
	\end{align*}
	we define the map
	\begin{align*}
		F_1:\qquad T^*(\Hom(\C^2,\C^4))\ &\longrightarrow\ \Hom(\C^2,\C^4\oplus(\C^4)^*)\\
		(v_1\oplus v_2)\oplus(v_1^*\oplus v_2^*)\ &\longmapsto\ (v_1\oplus (-v_2^*))\oplus(v_2\oplus v_1^*).
	\end{align*}
	And we check $F_1$ is a symplectic isomorphism:
	\begin{align*}
		&\ \ \ F_1^*\omega_1\big((v_1\oplus v_2)\oplus(v_1^*\oplus v_2^*),(w_1\oplus w_2)\oplus(w_1^*\oplus w_2^*)\big)\\
		& = \omega_1\big((v_1\oplus (-v_2^*))\oplus(v_2\oplus v_1^*),(w_1\oplus (-w_2^*))\oplus(w_2\oplus w_1^*)\big)\\
		& =\big(\!-\!v_2^*(w_2)+w_1^*(v_1)\big)-\big(v_1^*(w_1)-w_2^*(v_2)\big)\\
		& =(w_1^*\oplus w_2^*)(v_1\oplus v_2)-(v_1^*\oplus v_2^*)(w_1\oplus w_2),
	\end{align*}
	which is the natural symplectic form on $T^*(\Hom(\C^2,\C^4))$.
	Notice that all three spaces $$T^*V,\quad T^*(\Hom(\C^2,\C^4)),\quad \Hom(\C^2,\C^4\oplus(\C^4)^*)$$ have natural Hamiltonian $\SL_2$-actions. Since both $F_0$ and $F_1$ are $\SL_2$-equivariant, $F$ is also an $\SL_2$-equivariant symplectic isomorphism as
	$$
		F=F_1\circ F_0.\qedhere
	$$
\end{proof}
\begin{Corollary}\label{TO}
    The map $F$ induces an isomorphism between $\overline{T^*(\SL_3/U)}$ and $\Ominbar$ as affine varieties.
\end{Corollary}
\begin{proof}
This is a direct corollary of Propositions \ref{ominhamil} and \ref{Fthm}.
\end{proof}

\section{The Triality Action on $\so_8$}\label{triality}
The triality action was first discovered by E.\ Cartan in his 1925 paper \cite{Cartan} in which he constructed a certain $S_3$-subgroup of the automorphism group of $\so_8$, such that an order three element $Aut(\so_8)$ is constructed from the following matrix
$$
\begingroup
\renewcommand*{\arraystretch}{1.5}
\begin{pmatrix*}[r]
    -\frac{1}{2} && -\frac{1}{2} && -\frac{1}{2} && -\frac{1}{2}\\
    \frac{1}{2} && \frac{1}{2} && -\frac{1}{2} && -\frac{1}{2}\\
    \frac{1}{2} && -\frac{1}{2} && \frac{1}{2} && -\frac{1}{2}\\
    \frac{1}{2} && -\frac{1}{2} && -\frac{1}{2} && \frac{1}{2}
\end{pmatrix*}.
\endgroup
$$
Fix the following set of simple roots for the root system of type $D_4$,
\begin{align*}
	\alpha_1 =\epsilon_1-\epsilon_2,\ 
	\alpha_2 =\epsilon_2-\epsilon_3,\ 
	\alpha_3 =\epsilon_3-\epsilon_4,\ 
	\alpha_4 =\epsilon_3+\epsilon_4.
\end{align*}
So that the Dynkin Diagram $D_4$ is labeled as follows.
    $$
     \dynkin[text style/.style={scale=1},label,label macro/.code={\alpha_{\drlap{#1}}},scale=1.2,edge
length=0.9cm]D4
    $$
Choose a Chevalley basis $$\{X_\alpha\}_{\alpha\in\Delta^+},\{Y_{-\alpha}\}_{\alpha\in\Delta^+},\{H_{\alpha_i}\}_{i=1,2,3,4}.$$ of $\so_8=\Wedge^2(\C^8)$ as follows
\begin{align*}
	& X_{\alpha_1}=e_1\wedge e_7 
	&& Y_{-\alpha_1}=e_2\wedge e_8\\
	& X_{\alpha_2}=e_2\wedge e_6
	&& Y_{-\alpha_2}=e_3\wedge e_7\\
	& X_{\alpha_3}=e_3\wedge e_5
	&& Y_{-\alpha_3}=e_4\wedge e_6\\
	& X_{\alpha_4}=e_3\wedge e_4
	&& Y_{-\alpha_4}=e_5\wedge e_6\\
	& X_{\alpha_1+\alpha_2}=e_1\wedge e_6 
	&& Y_{-\alpha_1-\alpha_2}=e_3\wedge e_8\\
	& X_{\alpha_2+\alpha_3}=e_2\wedge e_5 
	&& Y_{-\alpha_2-\alpha_3}=e_4\wedge e_7\\
	& X_{\alpha_2+\alpha_4}=e_2\wedge e_4
	&& Y_{-\alpha_2-\alpha_4}=e_5\wedge e_7\\
	& X_{\alpha_2+\alpha_3+\alpha_4}=e_2\wedge e_3\hspace{-3em}
	&& Y_{-\alpha_2-\alpha_3-\alpha_4}=e_6\wedge e_7\\
	& X_{\alpha_1+\alpha_2+\alpha_4}=e_1\wedge e_4\hspace{-3em} 
	&& Y_{-\alpha_1-\alpha_2-\alpha_4}=e_5\wedge e_8\\
	& X_{\alpha_1+\alpha_2+\alpha_3}=e_1\wedge e_5\hspace{-3em} 
	&& Y_{-\alpha_1-\alpha_2-\alpha_3}=e_4\wedge e_8\\
	& X_{\alpha_1+\alpha_2+\alpha_3+\alpha_4}=e_1\wedge e_3\hspace{-3em}
	&& Y_{-\alpha_1-\alpha_2-\alpha_3-\alpha_4}=e_6\wedge e_8\\
	& X_{\alpha_1+2\alpha_2+\alpha_3+\alpha_4}=e_1\wedge e_2\hspace{-3em}
	&& Y_{-\alpha_1-2\alpha_2-\alpha_3-\alpha_4}=e_7\wedge e_8\\
	& H_{\alpha_1}=e_1\wedge e_8-e_2\wedge e_7\hspace{-3em}
	&& H_{\alpha_3}=e_3\wedge e_6-e_4\wedge e_5\\
    &  H_{\alpha_2}=e_2\wedge e_7-e_3\wedge e_6\hspace{-3em}
    && H_{\alpha_4}=e_3\wedge e_6+e_4\wedge e_5
\end{align*}    %<-------------------

Since the Dynkin diagram $D_4$ 
has an $S_3$-symmetry and each automorphism of the Dynkin diagram does lift uniquely to a Lie algebra automorphism (c.f. Theorem 2.108 in \cite{Knapp}, and Lemma 2.6 in \cite{EL}), the existence of the triality action might seem easy at first glance. But, a priori, it is not clear that the these lifts gives an $S_3\hookrightarrow\mathrm{Aut}(\so_8)$. For example, not every lift of the cyclic permutation of the simple roots $\alpha_1,\alpha_3,\alpha_4$ has order three in $\mathrm{Aut}(\so_8)$.
In the work of \cite{MW}, a lifted triality action on the Lie algebra $\so(4,4)$ is constructed, but the authors proved their result by an explicit calculation in coordinates by computer.
In this section, a new interpretation (see Definition \ref{act}) of the triality action on $\so_8=\so(4,4)_\C$ is given by applying the decomposation of $\so_8$ into irreducible $\slc$-representations.

\begin{Lemma}
    The map $\eta$ defined in \eqref{eta} induces an isomorphism 
$$\Wedge^2\eta:\Wedge^2(\C^3\oplus\C\oplus\C^*\oplus(\C^3)^*)\longrightarrow\Wedge^2(\C^8)=\so_8$$
of $\slc$-representations, where $\so_8$ is viewed as an $\slc$-representation with respect to the adjoint action under the embedding 
$$
    \varphi_1:\slc\longhookrightarrow\so_8.
$$
which restricts to an embedding $\varphi_1|_{\h_{\slc}}:\Lieh_{\slc}\hookrightarrow\Lieh_{\so_8}$ such that
$$
\varphi_1^*(\alpha_2)=\alpha,\textrm{ and } \varphi_1^*(\alpha_1+\alpha_2+\alpha_3+\alpha_4)=\beta.
$$
\end{Lemma}
\begin{proof}
    First, let $V_0\coloneqq\C^3\oplus\C\oplus\C^*\oplus(\C^3)^*$, and $<\,,\,>$ denote the inner product on $V_0$ given by the natural pairing in between $\C^3\oplus\C$ and $\C^*\oplus(\C^3)^*$. Then
    $$
        \Wedge^2(\C^3\oplus\C\oplus\C^*\oplus(\C^3)^*)=\so(V_0,<\,,\,>).
    $$
    Since $\eta$ is an isometry between $(V_0,<\,,\,>)$ and $(\C^8,(\,,\,))$, the map $\Wedge^2\eta$ is a Lie algebra isomorphism between $\so(V_0,<\,,\,>)$ and $\so_8$.
    We identify $$\slc=\{A\in\C^3\otimes(\C^3)^*\suchthat \tr(A)=0\}\subset\Wedge^2(\C^3\oplus\C\oplus\C^*\oplus(\C^3)^*).$$
    Then $\slc$ becomes an Lie subalgebra of $\so(V_0,<\,,\,>)$, and the $\slc$-representation structure on $\Wedge^2(\C^3\oplus\C\oplus\C^*\oplus(\C^3)^*)$ is given by the adjoint action.
    Then we have the restriction of $\Wedge^2\eta$ to $\slc$ equals to the embedding $\varphi_1$. For example, for positive root vectors,
    \begin{align*}
        \Wedge^2\eta\,(e_1\wedge e^*_2)& =e_2\wedge e_6= X_{\alpha_2},\\
        \Wedge^2\eta\,(e_2\wedge e^*_3)& =e_3\wedge (-e_1)= X_{\alpha_1+\alpha_2+\alpha_3+\alpha_4},\\
        \Wedge^2\eta\,(e_1\wedge e^*_3)& =e_2\wedge (-e_1)= X_{\alpha_1+2\alpha_2+\alpha_3+\alpha_4}.
    \end{align*}
    Since the $\slc$-representation structure on $\so_8$ is given by the Lie algebra embedding $\varphi_1$, our statement follows. 
\end{proof}
Decompose $\so_8=\Wedge^2(\C^3\oplus\C\oplus\C^*\oplus(\C^3)^*)$ into irreducible $\slc$-representations, 
\begin{equation}\label{decomp}
\so_8=\slc\oplus\C_{\textrm{trace}}\oplus\C\oplus\C^3\oplus\C^3\oplus\C^3\oplus(\C^3)^*\oplus(\C^3)^*\oplus(\C^3)^*.
\end{equation}
Let $\h=\{(c_1,c_2,c_3)\in\C^3\,|\,c_1+c_2+c_3=0\}$. Define $\varphi_2:\h\rightarrow\so_8$ by 
\begin{align*}
	\varphi_2(c_1,c_2,c_3)\coloneqq\ & c_1H_{\alpha_1}+c_2H_{\alpha_2}+c_3H_{\alpha_3}\\
	=\ & c_1(H_1-H_2)+c_2(H_3-H_4)+c_3(H_3+H_4)\\
	=\ & -c_1(H_2+H_3-H_1)+(c_3-c_2)H_4\in\C_{\textrm{trace}}\oplus\C,
\end{align*}
where $H_i\coloneqq e_i\wedge e_{9-i}$ for each $i$.
So the image of $\varphi_2$ is precisely the sub-representation $\C_{\textrm{trace}}\oplus\C\subset\so_8.$
Let $V_1,V_2,V_3$ (resp.\ $V^*_1,V^*_2,V^*_3$) be the three copies of $\C^3$ (resp.\ $(\C^3)^*$) in the decomposition (\ref{decomp}) whose highest weight vectors are root vectors for $\so_8$ corresponding to the roots
$\alpha_1+\alpha_2, \alpha_2+\alpha_3,\alpha_2+\alpha_4$ (resp.\ $\alpha_2+\alpha_3+\alpha_4, \alpha_1+\alpha_2+\alpha_4, \alpha_1+\alpha_2+\alpha_3$). Identify
 $V_1,V_2,V_3$ (\textrm{resp.\ }$V^*_1,V^*_2,V^*_3$) via the unique $\slc$-equivariant linear isomorphisms which map the above choice of highest weight vectors to each other.
 Explicitly,
\begin{align*}
& V_1=\C\langle -X_{\alpha_1+\alpha_2},X_{\alpha_1},Y_{-\alpha_2-\alpha_3-\alpha_4}\rangle=\C^3,\\
& V_2=\C\langle X_{\alpha_2+\alpha_3},X_{\alpha_3},Y_{-\alpha_1-\alpha_2-\alpha_4}\rangle=\C^3,\\
& V_3=\C\langle X_{\alpha_2+\alpha_4},X_{\alpha_4},Y_{-\alpha_1-\alpha_2-\alpha_3}\rangle=\C^3,\\
& V_1^*=\C\langle X_{\alpha_2+\alpha_3+\alpha_4},Y_{-\alpha_1},-Y_{-\alpha_1-\alpha_2}\rangle=(\C^3)^*,\\
& V_2^*=\C\langle X_{\alpha_1+\alpha_2+\alpha_4},Y_{-\alpha_3},Y_{-\alpha_2-\alpha_3}\rangle=(\C^3)^*,\\
& V_3^*=\C\langle X_{\alpha_1+\alpha_2+\alpha_3},Y_{-\alpha_4},Y_{-\alpha_2-\alpha_4}\rangle=(\C^3)^*.
\end{align*}
So we have an isomorphism of $\slc$-representations
\begin{equation}\label{phi}
       \varphi: \slc\oplus\h\oplus V_1\oplus V_2\oplus V_3\oplus V^*_1\oplus V^*_2\oplus V^*_3\longrightarrow\so_8
\end{equation}

\begin{Definition}\label{act}
    We define an $S_3$-action on $\slc\oplus\h\oplus V_1\oplus V_2\oplus V_3\oplus V^*_1\oplus V^*_2\oplus V^*_3$ by fixing $\slc$, acting on $\h$ as the Weyl group $S_3$-action, and permuting subscripts of $V_i$ and $V^*_i$. Then under the identification $\varphi$, we have an $S_3$-action $\act$ on $\so_8$.
\end{Definition}
\begin{Theorem}
    The $S_3$-action $\act$ gives an embedding $S_3\hookrightarrow \mathrm{Aut}(\so_8)$ which coincides with the triality action.
\end{Theorem}

\begin{proof}
	Let $(M,{c},{u},{u}^*)$ be an element in $\slc\oplus\h\oplus(\C^3)^3\oplus((\C^3)^*)^3$. Express each components in terms of the Chevalley basis, we see the $S_3$-action $\act$ does permutes the root vectors accordingly (i.e. fixing $\alpha_2$, and permutes $\alpha_1,\alpha_3,\alpha_4$). So it suffices to show that this action preserves the Lie bracket of $\so_8$.
There is a linear action of $(M,{c},{u},{u}^*)$ on elements $(v,a,b,v^*)\in\C^3\oplus\C\oplus\C^*\oplus(\C^3)^*$ coming from the linear action of $\so_8$ on $\C^8$ given by
	$$
	(M,c,u,u^*).\begin{pmatrix}
		v\\
		a\\
		b\\
		\ v^*
	\end{pmatrix}=
	\begin{pmatrix}
		M v+c_1v+au_2+bu_3+v^*\wedge u_1^*\\
		u_2^*(v)+(c_3-c_2)a-v^*(u_3)\\
		u_3^*(v)-(c_3-c_2)b-v^*(u_2)\\
		-v^* M-c_1v^*-au_3^*-bu_2^*+v\wedge u_1
	\end{pmatrix},
	$$
which can be applied to compute the Lie bracket structure on $\slc\oplus\h\oplus(\C^3)^3\oplus((\C^3)^*)^3$.

	Let $A=(M_A,{c}_A,{u}_A,{u}^*_A), B=(M_B,c_B,u_B,u^*_B)\in\slc\oplus\h\oplus(\C^3)^3\oplus((\C^*)^3)^3.$ Then the Lie bracket can be written as $[A,B]=(M_C,c_C,u_C,u^*_C)$ where
	\renewcommand{\arraystretch}{1.5}
	\begin{align*}
		\!\!M_C& =[M_A,M_B]+\sum_{i=1}^{3}(u_{A,i}\otimes u^*_{B,i}-u_{B,i}\otimes u^*_{A,i})+\frac{1}{3}\sum_{i=1}^{3}(u^*_{A,i}(u_{B,i})\!-\!u^*_{B,i}(u_{A,i}))\,\Id_{\C^3}\\
		c_C& =
		\begin{pmatrix}
			-2/3 && 1/3 && 1/3\\
			1/3 && -2/3 && 1/3\\
			1/3 && 1/3 && -2/3
		\end{pmatrix}
		\begin{pmatrix}
			u^*_{A,1}(u_{B,1})-u^*_{B,1}(u_{A,1}) \\
			u^*_{A,2}(u_{B,2})-u^*_{B,2}(u_{A,2}) \\
			u^*_{A,3}(u_{B,3})-u^*_{B,3}(u_{A,3})
		\end{pmatrix}\\
		u_C& =		\begin{pmatrix}
			M_A\,u_{A,1}-M_B\,u_{B,1}+u^*_{A,2}\wedge u^*_{B,3}-u^*_{B,2}\wedge u^*_{A,3}+2c_{A,1}\,u_{B,1}-2c_{B,1}\,u_{A,1}\\
			M_A\,u_{A,2}-M_B\,u_{B,2}+u^*_{A,3}\wedge u^*_{B,1}-u^*_{B,3}\wedge u^*_{A,1}+2c_{A,2}\,u_{B,2}-2c_{B,2}\,u_{A,2}\\
			M_A\,u_{A,3}-M_B\,u_{B,3}+u^*_{A,1}\wedge u^*_{B,2}-u^*_{B,1}\wedge u^*_{A,2}+2c_{A,3}\,u_{B,3}-2c_{B,3}\,u_{A,3}\\
		\end{pmatrix}\\
		u^*_C& =		-\begin{pmatrix}
			u^*_{A,1}M_A-u^*_{B,1}M_A+u_{A,2}\wedge u_{B,3}-u_{B,2}\wedge u_{A,3}+2c_{A,1}\,u^*_{B,1}-2c_{B,1}\,u^*_{A,1}\\
			u^*_{A,2}M_A-u^*_{B,2}M_A+u_{A,3}\wedge u_{B,1}-u_{B,3}\wedge u_{A,1}+2c_{A,2}\,u^*_{B,2}-2c_{B,2}\,u^*_{A,2}\\
			u^*_{A,3}M_A-u^*_{B,3}M_A+u_{A,1}\wedge u_{B,2}-u_{B,1}\wedge u_{A,2}+2c_{A,3}\,u^*_{B,3}-2c_{B,3}\,u^*_{A,3}\\
		\end{pmatrix}
	\end{align*}
	Then we see that the Lie bracket is manifestly equivariant under the lifted triality action. 
\end{proof}
In fact, the embedding of $\Tslc$ into the LHS of \eqref{phi} is a special case of the following conjecture due to Ginzburg and Kazhdan.
\begin{Conjecture}
        There is an $W$-equivariant embedding of the affine variety
$$
    \overline{T^*(G/U)}\hookrightarrow \Lieg\oplus\Lieh\oplus\bigoplus_{\varpi\in\textrm{Fund.\ Weights}}V_{\varpi}^{\,|\textrm{the }W\textrm{-orbit of }\varpi|},
$$
such that it restricts to the embedding 
$\overline{G/U}\hookrightarrow\bigoplus_{\varpi\in\textrm{Fund.\ Weights}}V_{\varpi}$, where the Weyl group $W$ acts on $\overline{T^*(G/U)}$ via the Gelfand--Graev action, acts on $\Lieh$ as the Weyl group, and permutes the copies of the fundamental representation $V_{\varpi}$.
\end{Conjecture}

\begin{Remark}
Fix the standard basis $e_1,e_2,e_3$ of $\C^3$, and identify
$$
\Wedge^3\C^3=\Wedge^3(\C^3)^*=\C,\ \Wedge^2\C^3=(\C^3)^*,\ \Wedge^2(\C^3)^*=\C^3.
$$
For any given $A\in\End(\C^3)$ and $v\in\C^3$ we define $A\wedge v\in\Sym^2(\C^3)^*$ by
$$
    (A\wedge v)(w_1,w_2)=(Aw_1)\wedge v\wedge w_2+(Aw_2)\wedge v\wedge w_1.
$$
Similarly we can define $A\wedge v^*\in\Sym^2\C^3$ for $v^*\in(\C^3)^*$.
Let $(M,{c},{u},{u}^*)$ be an element in $\slc\oplus\h\oplus(\C^3)^3\oplus((\C^3)^*)^3$.
Then 
$\varphi(M,{c},{u},{u}^*)\in\Ominbar$ if and only if:
\begin{align*}
	& u_i^*(u_j)=0\textrm{ for all distinct $i,j\in\{1,2,3\}$},\\
	& u^*_1(u_1)=(c_1-c_2)(c_1-c_3),\textrm{ and its cyclic permutations (``c.p." in short)},\\
	& u_1\wedge u_2 = (c_1-c_2)u_3^*,\textrm{ and its c.p.},\\
	& u_1^*\wedge u_2^* = (c_1-c_2)u_3,\textrm{ and its c.p.},\\
	& (M-c_3\,\Id_{\C^3})\wedge u_3+u_1^*\cdot u_2^*=0, \textrm{ and its c.p.},\\
	& (M-c_3\,\Id_{\C^3})\wedge u_3^*+u_1\cdot u_2=0, \textrm{ and its c.p.},\\
	& (M-c_3\,\Id_{\C^3})^2+u_1\otimes u_1^*+u_2\otimes u_2^*-u_3\otimes u_3^*+u_3^*(u_3)\,\Id_{\C^3}=0, \textrm{ and its c.p.}.
\end{align*}
It would be interesting to find all the relations for $\overline{T^*(\mathrm{SL}_n/U)}$ if the conjecture holds.
\end{Remark}

\section{The Gelfand-Graev Action}\label{ggaction}

We first recall the reconstruction in \cite{Wang} of the Gelfand-Graev action in Corollary 1.3.4 of \cite{GK}.
	Let $B=(\alpha,\beta)\in N\coloneqq\mu_H^{-1}(0)$.
	Let $k\in\{1,2,3,\cdots,n-1\}$. Define
	\begin{align*}
		& \out_k(B) \coloneqq\alpha_k\oplus\beta_{k-1}\in\Hom(\C^k,\C^{k+1}\oplus\C^{k-1}),\\
		& \inn_k(B) \coloneqq\beta_k\oplus(-\alpha_{k-1})\in\Hom(\C^{k+1}\oplus\C^{k-1},\C^k).
	\end{align*}
Then
$$
	\inn_k(B)\,\out_k(B)=\lambda_k\,\textrm{Id}_{\C^k}.
$$
	We also define $Z_k$ to be a subvariety of $N\times N$ consisting of pairs $$(B,B')=((\alpha,\beta),(\alpha',\beta')),$$ such that\\
(1)\quad $\forall j\notin\{k-1,k\},\ \alpha_j=\alpha'_j,\textrm{ and }\beta_j=\beta'_j.$\\
(2)\quad The following short sequence is exact
\begin{equation*}
    \xymatrix@=4em{
{0} \ar[r]
& {\C^k} \ar[r]^{\out_k(B')\qquad}
& {\C^{k+1}\oplus\C^{k-1}}  \ar[r]^{\qquad\inn_k(B)}
& {\C^k}\ar[r]
& 0}.
\end{equation*}
Moreover, for each $j$ we fix some volume form $\vol_j\in\Wedge^j\C^j$, and require
$$
	\out_k(B')(\vol_k)\wedge\inn_k(B)^{-1}(\vol_k)=\vol_{k-1}\wedge \vol_{k+1}\in\Wedge^{2k}(\C^{k+1}\oplus\C^{k-1}).
$$
(Since the short sequence is exact, we may take any element of $\inn_k(B)^{-1}(\vol_k)$.)\\
(3)
$$
	\out_k(B')\,\inn_k(B')=\out_k(B)\,\inn_k(B)-\lambda_k\,\textrm{Id}_{\C^{k+1}\oplus\C^{k-1}}.
$$

\begin{Theorem}[\cite{Wang}]
For any simple reflection $s_k\in W$, we can construct an automorphism $S_k$ on $N\sslash H$, such that for all $(B,B')\in Z_k$,
	  $$
	  		S_k(\pi(B))=\pi(B').
	  $$
And the automorphisms $S_k$ satisfy the braid relations, hence we constructed a $W$-action on $N\sslash H$. This $W$-action coincides with the quasi-classical Gelfand-Graev action for $\Tsln$.
\end{Theorem}

\begin{Theorem}
    In the case of $n=3$, and under the identification $$\Tslcbar=\Ominbar\subset\so_8,$$ the Gelfand-Graev action coincides with the triality $S_3$-action $\act$ on $\so_8$ restricted on $\Ominbar$.
\end{Theorem}
\begin{proof}
	Since Lie algebra automorphism preserves the minimal orbit, so the triality action can be restricted to an $S_3$ action on $\Ominbar$. It suffices to check for two of the simple transpositions $s_1=(23),s_2=(13)\in S_3$, the triality action satifies the condition that for each $k\in\{1,2\}$, and $(B,B')\in Z_k$, we have
	$$
		S_k(\pi(B))=\pi(B').
	$$
	Since the Gelfand-Graev action is uniquely determined by its restriction on the regular semisimple open part, we take $(M,{c},{u},{u}^*)$ to be an element in $\Omin\subset\slc\oplus\h\oplus(\C^3)^3\oplus((\C^3)^*)^3$ with distinct $c_1, c_2, c_3$.
	Then the corresponding isotropic decomposable element in $\Wedge^2 {\left(\C^3\oplus\C\oplus\C^*\oplus(\C^3)^*\right)}$ is
	$$
		%	(v_1+v_1^*)\wedge(v_2+v_2^*) & = \frac{1}{a_1b_2-a_2b_1}(a_2(v_1+v_1^*)-a_1(v_2+v_2^*))\wedge(b_2(v_1+v_1^*)-b_1(v_2+v_2^*))\\
		\left(u_3\oplus 0 \oplus -(c_3-c_2)\oplus -u_2^*\right)\wedge\left(\frac{u_2}{c_3-c_2}\oplus 1 \oplus 0 \oplus -\frac{u_3^*}{c_3-c_2}\right).
	$$
	So with this choice of representative, we have a lift $B=(\alpha,\beta)\in N$:
$$
	\hspace{0em}\alpha_1=
	\begin{pmatrix}
		0\\
		\!c_3\!-\!c_2\!
	\end{pmatrix}\!,\
	\beta_1=
	\begin{pmatrix}
		0 & 1
	\end{pmatrix}\!,\qquad\quad\ 
 $$
 $$
	\alpha_2=
	\vcenter{\begin{tikzpicture}\hspace{-3em}
			\matrix (m)[matrix of math nodes, nodes in empty cells, left delimiter={(},right delimiter={)},minimum width=3.2em,minimum height=2.2em)]
			{
				\   & \    \\
				u_3 & \!\!\!\displaystyle\frac{u_2}{c_3\!-\!c_2} \!\!\!  \\
				\   & \     \\
			} ;
			\draw (m-1-1.north west) rectangle (m-3-2.south east);
			\draw (m-1-1.north east) -- (m-3-1.south east);	
	\end{tikzpicture}}\hspace{-27.25 em},\ 
	\beta_2=
\hspace{-3em}
\vcenter{\begin{tikzpicture}
		\matrix (m)[matrix of math nodes, nodes in empty cells, left delimiter={(},right delimiter={)},minimum width=3.6em,minimum height=2.75em)]
		{
			\quad\ \displaystyle\ \ -\frac{u_3^*}{c_3-c_2} \qquad \\
			\qquad\ \ \ \ u_2^* \qquad\ \ \ \\
		} ;
		\draw (m-2-1.south west) rectangle (m-1-1.north east);
		\draw (m-2-1.north west) -- (m-1-1.south east);	
\end{tikzpicture}}
\hspace{-22em}.
$$

First, for $k=1$, we have
	$$
	\out_1(B)=
		\begin{pmatrix}
			0\\
			c_3-c_2
		\end{pmatrix},\quad
	\inn_1(B)=
		\begin{pmatrix}
			0 & 1
		\end{pmatrix},\quad
	\lambda_1=c_3-c_2.
	$$
	Apply the action by $s_1=(23)$, and rewrite the bivector so that the first defining property of $Z_1$ holds.
	\begin{align*}
		&\ \ \ \left(u_2\oplus 0 \oplus -(c_2-c_3)\oplus -u_3^*\right)\wedge\left(\frac{u_3}{c_2-c_3}\oplus 1 \oplus 0 \oplus -\frac{u_2^*}{c_2-c_3}\right)\\
		& = \left(u_3\oplus (c_2-c_3) \oplus 0\oplus -u_2^*\right)\wedge\left(\frac{u_2}{c_3-c_2}\oplus 0 \oplus 1 \oplus -\frac{u_3^*}{c_3-c_2}\right).
	\end{align*}
So we have
$$
	\quad\ \out_1(B')=
	\begin{pmatrix}
		1\\
		0
	\end{pmatrix},\quad
	\inn_1(B')=\begin{pmatrix}
		\,c_2\!-\!c_3 & 0\,
	\end{pmatrix}.
$$
So $$
\out_1(B')\,\inn_1(B')=\out_1(B)\,\inn_1(B)-\lambda_1\,\Id_{\C^2}
.$$ So $(B,B')\in Z_1$, and the action of $s_1$ on $\Ominbar$ coincides with the action of $S_1$ on $\Tslcbar$.

Next, for $k=2$, we have $\lambda_2=c_1-c_3,$ and
$$
\out_2(B)=\hspace{-1.5em}
\vcenter{\begin{tikzpicture}
		\matrix (m)[matrix of math nodes, nodes in empty cells, left delimiter={(},right delimiter={)},minimum width=3.2em,minimum height=2em)]
		{
			\   & \    \\
			u_3 & \!\!\!\displaystyle\frac{u_2}{c_3\!-\!c_2} \!\!\!  \\
			\   & \     \\
			0 & \!\!\!1  \!\!\! \\
		} ;
		\draw (m-1-1.north west) rectangle (m-4-2.south east);
		\draw (m-4-1.north west) -- (m-4-2.north east);
		\draw (m-1-1.north east) -- (m-4-1.south east);	
\end{tikzpicture}}	
\hspace{-29.5em},\qquad \inn_2(B)=\hspace{-1.5em}
\vcenter{\begin{tikzpicture}
		\matrix (m)[matrix of math nodes, nodes in empty cells, left delimiter={(},right delimiter={)},minimum width=3.2em,minimum height=2.75em)]
		{
			\quad\ \  \displaystyle-\frac{u_3^*}{c_3-c_2} \qquad  & 0\\
			\qquad\qquad u_2^* \ \ \ \qquad\ \ \ & c_2\!-\!c_3\\
		} ;
		\draw (m-2-1.south west) rectangle (m-1-2.north east);
		\draw (m-1-2.north west) -- (m-2-2.south west);	
		\draw (m-2-1.north west) -- (m-2-2.north east);	
\end{tikzpicture}}\hspace{-23em}.
$$
Applying the triality action by element $s_2=(13)\in S_3$, we get $B''$ such that
$$
\out_2(B'')=\hspace{-1.5em}
\vcenter{\begin{tikzpicture}
		\matrix (m)[matrix of math nodes, nodes in empty cells, left delimiter={(},right delimiter={)},minimum width=3.2em,minimum height=2em)]
		{
			\   & \    \\
			u_1 & \!\!\!\displaystyle\frac{u_2}{c_1\!-\!c_2} \!\!\!  \\
			\   & \     \\
			0 & \!\!\!1  \!\!\! \\
		} ;
		\draw (m-1-1.north west) rectangle (m-4-2.south east);
		\draw (m-4-1.north west) -- (m-4-2.north east);
		\draw (m-1-1.north east) -- (m-4-1.south east);	
\end{tikzpicture}}	
\hspace{-29.5em},\qquad \inn_2(B'')=\hspace{-1.5em}
\vcenter{\begin{tikzpicture}
		\matrix (m)[matrix of math nodes, nodes in empty cells, left delimiter={(},right delimiter={)},minimum width=3.6em,minimum height=2.75em)]
		{
			\quad\ \  \displaystyle-\frac{u_1^*}{c_1-c_2} \qquad  & 0\\
			\qquad\qquad u_2^* \ \ \ \qquad\ \ \ & c_2\!-\!c_1\\
		} ;
		\draw (m-2-1.south west) rectangle (m-1-2.north east);
		\draw (m-1-2.north west) -- (m-2-2.south west);	
		\draw (m-2-1.north west) -- (m-2-2.north east);	
\end{tikzpicture}}
\hspace{-22.5em}.
$$
Then we can check
$$
\out_2(B)\,\inn_2(B) 
=\hspace{-1.5em}\vcenter{\begin{tikzpicture}
		\matrix (m)[matrix of math nodes, nodes in empty cells, left delimiter={(},right delimiter={)},minimum width=3.3em,minimum height=2.2em)]
		{
			\ \hspace{6em}\ & \    \\
			-M+c_1\Id_{\C^3} & u_2   \\
			\ \hspace{6em}\ & \     \\
			u^*_2 & c_2\!-\!c_3 \\
		} ;
		\draw (m-1-1.north west) rectangle (m-4-2.south east);
		\draw (m-3-1.south west) -- (m-3-2.south east);
		\draw (m-1-2.north west) -- (m-4-2.south west);	
\end{tikzpicture}}\hspace{-25em},
$$
$$
\out_2(B'')\,\inn_2(B'') =\hspace{-1.5em}\vcenter{\begin{tikzpicture}
		\matrix (m)[matrix of math nodes, nodes in empty cells, left delimiter={(},right delimiter={)},minimum width=3.3em,minimum height=2.2em)]
		{
			\ \hspace{6em}\ & \    \\
			-M+c_3\Id_{\C^3} & u_2   \\
			\ \hspace{6em}\ & \     \\
			u^*_2 & c_2\!-\!c_1 \\
		} ;
		\draw (m-1-1.north west) rectangle (m-4-2.south east);
		\draw (m-3-1.south west) -- (m-3-2.south east);
		\draw (m-1-2.north west) -- (m-4-2.south west);	
\end{tikzpicture}}\hspace{-25em},
$$
where we have used the relations:
$$
\frac{u_2\wedge u^*_2-u_3\wedge u_3^*}{c_3-c_2}=-M+c_1\Id_{\C^3},
\qquad
\frac{u_2\wedge u^*_2-u_1\wedge u_1^*}{c_1-c_2}=-M+c_3\Id_{\C^3}.
$$
So 
$$
\out_2(B'')\,\inn_2(B'')=\out_2(B)\,\inn_2(B)-\lambda_2\,\Id_{\C^4}
$$
So $(B,B'')\in Z_2$.
This shows that the Gelfand-Graev action on $\Tslcbar$ coincides with the restriction of the triality action on $\Ominbar$.
\end{proof}

\section{Symplectic form on the smooth locus}\label{symplectic}
Let us recall the following known results from \cite{CG}.
Let $\Onilp$ be a nilpotent adjoint orbit in $\Lieg$.
Since $\Onilp$ is a nilpotent orbit, there is a $\C^*$-action on it such that $\omega^\textrm{KKS}=\mathcal{L}_{\textsl{Eu}}\omega^\textrm{KKS},$ where $\textsl{Eu}$ is the Euler vector field on $\Onilp$. Define $$\lambda^\textrm{KKS}\coloneqq\iota_{\textsl{Eu}}\omega^\textrm{KKS}.$$
So we have $$d\lambda^\textrm{KKS}=d\iota_{\textsl{Eu}}\omega^\textrm{KKS}=\mathcal{L}_{\textsl{Eu}}\omega^\textrm{KKS}=\omega^\textrm{KKS}.$$
\begin{Proposition}
    Let $X\in\Onilp$ and $Y\in\Lieg$.  Then
    $$\lambda^\textrm{KKS}_X(\xi_Y)=\kappa(X,Y).$$
\end{Proposition}
\begin{proof}
Since $X$ is nilpotent, we have $Y\in\Lieg_X$ if and only if $\kappa(X,Y)=0$. So the RHS of the formula only depends on $\xi_Y$ (i.e. independent of the choice of $Y$). Let $E\in\Onilp$. Fix some $\slb$-triple $(E,H,F)$ in $\Lieg$. Then we have
$$
    (\xi_H)_E=2\,\textsl{Eu}_{E}\in T_E\Onilp.
$$
Since $\lambda^\textrm{KKS}$ is invariant under the $G$-action, it suffices to prove the statement for $X=E$. Let $Y\in\Lieg$. We compute
$$
(\iota_{\textsl{Eu}}\omega_E^\textrm{KKS})(\xi_Y)=\frac{1}{2}\omega_E^\textrm{KKS}(\xi_Y,\xi_H)=\frac{1}{2}\kappa(E,[Y,H])=\frac{1}{2}\kappa([H,E],Y)=\kappa(E,Y).\qedhere
$$
\end{proof}

\begin{Theorem}
	The symplectic form on $\overline{T^*(\SL_3/U)}_\textrm{sm}$ coincides (up to a scalar multiple) with $\omega^\textrm{KKS}$ on $\Omin\subset\so_8$.
\end{Theorem}
\begin{proof}
	Since the symplectic form on $\overline{T^*(\SL_3/U)}_\textrm{sm}$ constructed in \textrm{Proposition \ref{symp}} is the extension of the symplectic form given by the Hamiltonian reduction of $N_\textrm{inj}$ by $H$, it coincides with the pull-back of the symplectic form on $(N_1\sslash \SL_2)_\textrm{sm}$ by the restriction of the isomorphism in \textrm{Corollary\ \ref{TO}} to the smooth points: $$F:\left(\Tslcbar\right)_\textrm{sm}\rightarrow \Omin.$$
	So it suffices to show that the symplectic form on $(N_1\sslash \SL_2)_\textrm{sm}$ coming from the Hamiltonian reduction coincides (up to a scalar multiple) with $\omega^\textrm{KKS}$. It suffices to check for one forms.
	
	Let $v_1\wedge v_2\in \Omin\subset\Wedge^2 {\C^8}$.
	Recall the element $\varphi_{v_1\wedge v_2}\in \so_8$ is defined in (\ref{varphi}) by the formula
	$$
	\varphi_{v_1\wedge v_2}(u)=(v_1,u)v_2-(v_2,u)v_1.
	$$
	Observe that the tangent space $T_{v_1\wedge v_2}\Omin$ is spanned by vectors given by infinitesimal actions of some $w_1\wedge w_2\in\so_8$. Let 
	$$
	\xi_{w_1\wedge w_2}=[w_1\wedge w_2,v_1\wedge v_2]
	$$
	be such a tangent vector.
	
	To calculate $\xi_{w_1\wedge w_2}$ we first we compute
	\begin{align*}
		(\varphi_{v_1\wedge v_2}\circ \varphi_{w_1\wedge w_2})(u)\ \ =\ \ \ &[(w_2,u)(v_2,w_1)-\,(w_1,u)(v_2,w_2)] v_1\\
		+&[(w_1,u)(v_1,w_2)-\,(w_2,u)(v_1,w_1)] v_2,\\
		(\varphi_{w_1\wedge w_2}\circ \varphi_{v_1\wedge v_2})(u)\ \ =\ \ \ &[(v_2,u)(w_2,v_1)-\,(v_1,u)(w_2,v_2)] w_1\\
		+&[(v_1,u)(w_1,v_2)-\,(v_2,u)(w_1,v_1)] w_2.
	\end{align*}
	Then we compute the Lie bracket $[\varphi_{v_1\wedge v_2}, \varphi_{w_1\wedge w_2}]$ and get
	\begin{align*}
		\xi_{w_1\wedge w_2}& =(v_2,w_2)\, w_1\wedge v_1-(v_1,w_2)\, w_1\wedge v_2-(v_2,w_1)\, w_2\wedge v_1+(v_1,w_1)\, w_2\wedge v_2\\
		& =\varphi_{w_1\wedge w_2}(v_1)\wedge v_2+ v_1\wedge \varphi_{w_1\wedge w_2}(v_2).
	\end{align*}
	The value of the one form $\lambda^\textrm{KKS}$ on $\Omin$ at $v_1\wedge v_2$ is
	\begin{align*}
		\lambda^\textrm{KKS}_{v_1\wedge v_2}(\xi_{w_1\wedge w_2})=\Tr(\varphi_{v_1\wedge v_2}\circ \varphi_{w_1\wedge w_2})=2((v_2,w_1)(v_1,w_2)-(v_1,w_1)(v_2,w_2)).
	\end{align*}

	Recall from (\ref{o1}), we have $\omega_1=(1/2)d\lambda'$, where $\lambda'$ is the $\SL_2$-invariant one form on $\Hom(\C^2,\C^8)=\C^8\oplus\C^8$ defined by
	$$\lambda'_{v_1\oplus v_2}(x_1\oplus x_2)=(v_1,x_2)-(v_2,x_1).$$

	Lift the tangent vector $\xi_{w_1\wedge w_2}$ to $(x_1\oplus x_2)\in T_{v_1\oplus v_2}(\C^8\oplus\C^8)$, where
	\begin{align*}
		x_1& =\varphi_{w_1\wedge w_2}(v_1)=(w_1,v_1)w_2-(w_2,v_1)w_1,\\
		x_2& =\varphi_{w_1\wedge w_2}(v_2)=(w_1,v_2)w_2-(w_2,v_2)w_1.
	\end{align*}
	So we have $$\lambda'_{v_1\oplus v_2}(x_1\oplus x_2)=2((w_1,v_2)(v_1,w_2)-(w_2,v_2)(v_1,w_1))=\lambda^\textrm{KKS}_{v_1\wedge v_2}(\xi_{w_1\wedge w_2}).\qedhere$$
\end{proof}

% Your bilbigraphy           %<-------------------

\printbibliography

\bigskip
\footnotesize
 \noindent Boming Jia, \textit{Email}: \texttt{jiabm@tsinghua.edu.cn}\\
\textsc{Yau Mathematical Sciences Center,\\
 Jingzhai 301, Tsinghua University,\\
 Beijing, 100084, China.}\par\nopagebreak
 
\end{document}